\documentclass[draft]{amsart}

\newtheorem{thm}{Theorem}
\newtheorem{lem}[thm]{Lemma}
\newtheorem{dfn}[thm]{Definition}

\theoremstyle{remark}
\newtheorem{rmk}[thm]{Remark}

\newcommand{\pv}{\mathbf{P}}

\numberwithin{thm}{section}

\begin{document}


\title{The Fatou Set for Critically Finite Maps}

\author{Feng Rong}

\address{Department of Mathematics, University of Michigan, Ann Arbor, MI 48109, USA}
\email{frong@umich.edu}

\begin{abstract}
It is a classical result in complex dynamics of one variable that the Fatou set for a critically finite map on $\pv^1$ consists of only basins of attraction for superattracting periodic points. In this paper we deal with critically finite maps on $\pv^k$. We show that the Fatou set for a critically finite map on $\pv^2$ consists of only basins of attraction for superattracting periodic points. We also show that the Fatou set for a $k-$critically finite map on $\pv^k$ is empty.
\end{abstract}

\maketitle

\section{Introduction}

A holomorphic map $f:\pv^k\rightarrow \pv^k$ is said to be critically finite if every component of the critical set for $f$ is periodic or preperiodic. In \cite{T:Comb}, Thurston has given a topological classification of critically finite maps on $\pv^1$. And it is well known that the Fatou set for a critically finite map on $\pv^1$ consists of only basins of attraction for superattracting periodic points, i.e. points $p$ with $f^n(p)=p$ and $(f^n)^\prime(p)=0$ for some $n\in \mathbf{N}$ (see \cite{M:One}). In this paper, we show that the same is also true for critically finite maps on $\pv^2$. More precisely, we have the following

\begin{thm}\label{T:Main1}
If $f:\pv^2\rightarrow \pv^2$ is a critically finite holomorphic map, then the Fatou set for $f$ consists of only basins of attraction for superattracting periodic points.
\end{thm}

With some extra assumptions, the above result has been obtained by Fornaess and Sibony (\cite{FS:DynamicsI}).

We will also study critically finite maps on $\pv^k$. In particular, we obtain the following (see Section \ref{S:Fatou} for precise definitions).

\begin{thm}\label{T:Main2}
Let $f:\pv^k\rightarrow \pv^k$ be a holomorphic map. If $f$ is $k-$critically finite, then the Fatou set for $f$ is empty.
\end{thm}

For $1-$critically finite maps on $\pv^1$ and $2-$critically finite maps on $\pv^2$, this was proved by Thurston (\cite{T:Comb}) and Ueda (\cite{U:Critical}), respectively. 

The author would like to thank John Erik Fornaess for his advice and encouragement.

\section{The Fatou Set for Critically Finite Maps}\label{S:Fatou}

Let $f:\pv^k\rightarrow \pv^k$ be a holomorphic map of (algebraic) degree $d > 1$.

Let $C_1$ be the critical set of $f$ given by
$$C_1 = \{ p\in \pv^k | rank(df(p)) < k\},$$
where $df(p)$ denotes the differential of $f$ at $p$.

We define the post-critical set $D_1$ of $f$ by
$$D_1 = \bigcup\limits_{j=1}^{\infty} f^j(C_1),$$
and the $\omega-$limit set $E_1$ of $f$ by
$$E_1 = \bigcap\limits_{j=1}^{\infty} f^j(\overline{D_1}).$$

By definition, a holomorphic map $f$ on $\pv^k$ is critically finite if the post-critical set $D_1$ is an analytic (hence algebraic) set in $\pv^k$. This is equivalent to saying that there is an integer $l\ge 1$ such that $D_1=\cup_{j=1}^l f^j(C_1)$. Hence, in the critically finite case, the set $D_1$ is an algebraic set of pure codimension 1.

Let us take a closer look at the structure of the post-critical set $D_1$ and the $\omega-$limit set $E_1$. If $f$ is critically finite, then $f^{j-1}(D_1) = \cup_{l=j}^{\infty}f^l(C_1),\ j=1,2,\cdots$, is a descending sequence of algebraic sets. Hence there is an integer $l_1\ge 1$ such that $f^{l_1-1}(D_1) = f^{l_1}(D_1) =\cdots$. Consequently $E_1 = f^{l_1-1}(D_1)$ is an algebraic set of pure codimension 1. We can decompose $E_1$ into $E_1^\prime\cup F_1$, where $F_1$ consists of those components in a critical cycle. (A periodic component $L$ is said to be in a \textit{critical cycle} if at least one of the forward images of $L$ under $f$ is contained in the critical set for $f$.)

\begin{dfn}
Let $f:\pv^k\rightarrow \pv^k$ be a holomorphic map. The map $f$ is said to be critically finite of order 1 if $D_1$, hence $E_1$, is algebraic. And $f$ is said to be 1-critically finite if $C_1$ and $E_1$ have no common irreducible component, i.e. $F_1=\emptyset$.
\end{dfn}

We can now make the following inductive definition (c.f. \cite{J:Critical}).

\begin{dfn}\label{D:Post}
Let $f:\pv^k\rightarrow \pv^k$ be a holomorphic map. Suppose $f$ is critically finite of order $n-1$, $1<n\le k$. Denote $C_n = C_1\cap E_{n-1}$, $D_n = \cup_{j=1}^{\infty}f^j(C_n)$, and $E_n = \cap_{j=1}^{\infty}f^j(\overline{D_n})$. We say that $f$ is critically finite of order $n$ if $D_n$, hence $E_n$, is algebraic. Let $l_n$ be the least integer such that $E_n = f^{l_n-1}(D_n) = \cup_{j=l_n}^{\infty}f^j(C_n)$. We can decompose $E_n$ into $E_n^\prime\cup F_n$, where $F_n$ consists of those components, of codimension less or equal to $n$, in a critical cycle. If in addition $f$ is $(n-1)-$critically finite, then we say that $f$ is $n-$critically finite if $E_n$ has no irreducible component contained in $C_1$, i.e. $F_n=\emptyset$.
\end{dfn}

Before we go further, let us recall some definitions and results from \cite{U:Critical}.

\begin{dfn}
Let $f:\pv^k\rightarrow \pv^k$ be a holomorphic map and let $U$ be a Fatou component, i.e. a connected component of the Fatou set for $f$. A holomorphic map 
$\varphi:U\rightarrow \pv^k$ is called a limit map on $U$ if there is a sequence $\{f^{n_j}|U\}$ which converges to $\varphi$ uniformly on compact sets in $U$. 
A point $q\in \pv^k$ is called a \textit{Fatou limit point} if there is a limit map $\varphi$ on a Fatou component $U$ such that $q\in \varphi(U)$. The set of all Fatou limit points is called the \textit{Fatou limit set}.
\end{dfn}

\begin{dfn}
A Fatou component $U$ is called a \textit{rotation domain} if the identity map $id_U:U\rightarrow U$ is a limit map on $U$.
\end{dfn}

\begin{dfn}\label{D:Ramification}
A point $q\in \pv^k$ is said to be a \textit{point of bounded ramification} with respect to $f$ if the following conditions are satisfied:\\
$(i)$ There is a neighborhood $W$ of $q$ such that $D_1\cap W$ is an analytic subset of $W$;\\
$(ii)$ There exists an integer $l$ such that, for every integer $j>0$ and every $p\in f^{-j}(q)$, the cardinality $\sharp(I)$ of the set 
$$I = \{i | 0\le i\le j-1, f^i(p)\in C_1\}$$
is not greater than $l$.
\end{dfn}

The following two theorems by Ueda are crucial. 

\begin{thm}\cite[Theorem 4.8]{U:Critical}\label{T:Rotation}
Suppose that $q\in \pv^k$ is a point of bounded ramification and also a Fatou limit point. Then $q$ is contained in a rotation domain.
\end{thm}

\begin{thm}\cite[Proposition 5.1, (1)]{U:Critical}\label{P:Rotation}
If $f:\pv^k\rightarrow \pv^k$ is critically finite, then there is no rotation domain.
\end{thm}

\begin{rmk}\label{R:Ramification}
Since the set $D_1$ is an analytic set in the critically finite case, condition $(i)$ in Definition \ref{D:Ramification} is automatically true. So we only need to check condition $(ii)$ in Definition \ref{D:Ramification} to see if a point $p\in \pv^k$ is of bounded ramification.
\end{rmk}

We need the following lemma,  whose proof is an elaboration of the proof of Lemma 5.7 in \cite{U:Critical}.

\begin{lem}
Let $f:\pv^k\rightarrow \pv^k$ be a holomorphic map. If $f$ is critically finite of order $n$, $1\le n < k$, then every point in $\pv^k\backslash E_n$ is a point of bounded ramification. If $f$ is critically finite of order $k$, then every point in $\pv^k\backslash F_k$ is a point of bounded ramification.
\end{lem}
\begin{proof}
First assume that $f$ is critically finite of order $n$, $1\le n < k$. Let $q\in \pv^k\backslash E_n$ and let $p\in f^{-j}(q)$ for some integer $j>0$. By Remark \ref{R:Ramification}, we only need to show that the cardinality $\sharp(I)$ of the set
$$I=\{i | 0\le i\le j-1, f^i(p)\in C_1\}$$
is not greater than some integer $l>0$. Let
$$I_m=\{i | 0\le i\le j-1, f^i(p)\in C_m\backslash C_{m+1}\},\ \ m=1,\cdots,n-1,$$
$$I_n=\{i | 0\le i\le j-1, f^i(p)\in C_n\}.$$
We claim that $\sharp(I_m)\le l_m$, $m=1,\cdots,n$.

For each $1\le m <n$, suppose that $I_m$ is non-empty and let $i_m$ be the least index in $I_m$. Then $f^{i_m}(p)\in C_m$. For $i\ge i_m+l_m$, we have $f^i(p)\in E_m$ and hence $f^i(p)\notin C_m\backslash C_{m+1}$. Thus $I_m$ is a subset of $\{i_m,\cdots,i_m+l_m-1\}$, and $\sharp(I_m)\le l_m$.

Now suppose that $I_n$ is non-empty and let $i_n$ be the least index in $I_n$. Then $f^{i_n}(p)\in C_n$. For $i\ge i_n+l_n$, we have $f^i(p)\in E_n$. Since $f^j(p)=q\notin E_n$, we have $i_n+l_n > j$. Thus $I_n$ is a subset of $\{i_n,\cdots,i_n+l_n-1\}$ and $\sharp(I_n)\le l_n$.

Next assume that $f$ is critically finite of order $k$. Let $q\in \pv^k\backslash F_k$ and let $p\in f^{-j}(q)$ for some integer $j>0$. Let
$$I_m=\{i | 0\le i\le j-1, f^i(p)\in C_m\backslash C_{m+1}\},\ \ m=1,\cdots,k-1,$$
$$I_k=\{i | 0\le i\le j-1, f^i(p)\in C_k\}.$$

For the same reason as above we have that $\sharp(I_m)\le l_m$ for $1\le m < k$. Now suppose that $I_k$ is non-empty and let $i_k$ be the least index in $I_k$. Then $f^{i_k}(p)\in C_k$. For $i\ge i_k+l_k$, we have $f^i(p)\in E_k$. Note that $f(F_k)=F_k$ and $(E_k\backslash F_k)\cap C_k =\emptyset$. Since $f^j(p)=q\notin F_k$, we have $i_k+l_k > j$. Thus $I_k$ is a subset of $\{i_k,\cdots,i_k+l_k-1\}$ and $\sharp(I_k)\le l_k$.
\end{proof}

Combining this lemma with Theorem \ref{T:Rotation} and \ref{P:Rotation}, we obtain the following

\begin{thm}\label{T:Limit}
Let $f:\pv^k\rightarrow \pv^k$ be a holomorphic map. If $f$ is critically finite of order $n$, $1\le n < k$, then the Fatou limit set is contained in $E_n$. If $f$ is critically finite of order $k$, then the Fatou limit set is contained in $F_k$.
\end{thm}

By definition, a $k-$critically finite map on $\pv^k$ has $F_k=\emptyset$. Therefore we obtain Theorem \ref{T:Main2} as a corollary to the above theorem.

Now let us turn our attention to critically finite maps on $\pv^2$. We say that a point $p$ is a superattracting periodic point for a homomorphic map $f$ on $\pv^2$ if there exists an $n\in \mathbf{N}$ such that $f^n(p)=p$ and both of the eigenvalues of the differential $df^n(p)$ are equal to zero, i.e. $df^n(p)$ is nilpotent. In \cite{FS:DynamicsI}, Forn$\ae$ss and Sibony obtained the following two theorems (adapted to our setting). (See \cite{FS:DynamicsI} for precise definitions.)

\begin{thm}\cite[Theorem 7.7]{FS:DynamicsI}\label{T:FS1}
Let $f:\pv^2\rightarrow \pv^2$ be a holomorphic map. Assume that $f$ is strictly critically finite and that $\pv^2\backslash \{C_1\cup D_1\}$ is hyperbolic. Then the only Fatou components of $f$ are basins of attraction for superattracting periodic points.
\end{thm}

\begin{thm}\cite[Theorem 7.8]{FS:DynamicsI}\label{T:FS2}
Let $f:\pv^2\rightarrow \pv^2$ be a holomorphic map. Assume that $f$ is strictly critically finite and that $F_1\neq \emptyset$. Then the only Fatou components of $f$ are basins of attraction for superattracting periodic points.
\end{thm}

Jonsson noted that a critically finite holomorphic map $f$ on $\pv^2$ is always strictly critically finite (\cite[Remark 2.10]{J:Critical}). We now prove Theorem \ref{T:Main1}.

\begin{proof}[Proof of Theorem \ref{T:Main1}]
If $f$ is not $1-$critically finite, we are done by Theorem \ref{T:FS2}. Therefore, we can assume that $f$ is $1-$critically finite. Then by \cite[Theorem 5.8]{U:Critical}, the Fatou limit set for $f$ consists of finitely many periodic critical points in $F_2$. Arguing as in the proof of Theorem \ref{T:FS1}, we are done.
\end{proof}

\begin{rmk}
Note that one does not necessarily have $df^n(p)=0$ for a superattracting periodic point $p$ of period $n$ for a critically finite map $f$ on $\pv^2$. For a simple example, let us look at the map $f:[z:w:t]\mapsto [z^2-wt:w^2:t^2]$. It is easy to check that $f$ is critically finite and $p=[0:0:1]$ is a superattracting fixed point for $f$. But $df(p)\neq 0$. Incidentally, maps similar to the example we just gave are of independent interest. Bonifant and Dabija (\cite[Theorem 4.1]{BD:Elliptic}) showed that an invariant critical component for a holomorphic map on $\pv^2$ must be a rational curve. While all known examples of critically finite maps on $\pv^2$ only have smooth rational curves as invariant critical components, we give here a family of critically finite maps on $\pv^2$ with singular rational curves as invariant critical components.

$$g_d:[z:w:t]\mapsto [z^d-w^{d-1}t:-w^d:-t^d],\ \ \ d > 2.$$

Note that $g_d$ maps the critical component $\{z=0\}$ to the singular rational curve $\{z^d=w^{d-1}t\}$ and maps $\{z^d=w^{d-1}t\}$ back to $\{z=0\}$. So $g_d^2$ will have $\{z^d=w^{d-1}t\}$ as a fixed critical component and obviously $g_d^2$ is a critically finite map.
\end{rmk}


\begin{thebibliography}{}

\bibitem[BD]{BD:Elliptic}
Bonifant, A., Dabija, M.;
\emph{Self-maps of $\pv^2$ with Invariant Elliptic Curves},
Contemp. Math., vol. 311 (2002), 1-25.

\bibitem[FS]{FS:DynamicsI}
Forn$\ae$ss, J.E., Sibony, S.;
\emph{Complex Dynamics in Higher Dimension. $I$},
Ast$\acute{e}$risque, 222 (1994), 201-231.

\bibitem[J]{J:Critical}
Jonsson, M.;
\emph{Some Properties of 2-Critically Finite Maps of $\pv^2$},
Ergordic Theory and Dynamical Systems, 18 (1998), 171-187.

\bibitem[M]{M:One}
Milnor, J.;
\emph{Dynamics in One Complex Variable},
Princeton Univ. Press, 3rd. ed., 2006. 

\bibitem[T]{T:Comb}
Thurston, W.;
\emph{On The Combinatorics and Dynamics of Rational Maps},
Preprint.

\bibitem[U]{U:Critical}
Ueda, T.;
\emph{Critical Orbits of Holomorphic Maps on Projective Spaces},
J. Geom. Anal., 8-2 (1998), 319-334.

\end{thebibliography}
\end{document}